\documentclass[11pt]{amsart}
\usepackage{amsmath,amssymb}
\usepackage[colorlinks=true,allcolors=blue]{hyperref}

\DeclareMathOperator{\Con}{Con}
\DeclareMathOperator{\cd}{cd}
\DeclareMathOperator{\con}{con}
\DeclareMathOperator{\Ji}{J}
\DeclareMathOperator{\qd}{qd}
\DeclareMathOperator{\Quo}{Quo}
\DeclareMathOperator{\Id}{Id}
\newcommand{\balpha}{\boldsymbol{\alpha}}
\newcommand{\bbeta}{\boldsymbol{\beta}}
\newcommand{\btheta}{\boldsymbol{\theta}}
\newcommand{\blambda}{\boldsymbol{\lambda}}
\newtheorem{theorem}{Theorem}[section]
\newtheorem{lemma}[theorem]{Lemma}
\newtheorem{corollary}[theorem]{Corollary}

\title{Congruence density under sublattice inclusion}
\author{George Gr\"atzer}
\date{}
\hypersetup{pdftitle={Congruence density under sublattice inclusion},
  pdfauthor={George Gratzer},pdfsubject={Finite lattice theory}}

\begin{document}
\begin{abstract}
Following G\'abor Cz\'edli, the \emph{congruence density} of a finite
nonempty lattice $A$ is $\cd(A)=|\Con A|/2^{|A|-1}$, where $\Con A$
is its congruence lattice and $|A|$ is its number of elements.
Cz\'edli proved that density cannot increase under several kinds of
sublattice extension and asked whether this holds for every sublattice
inclusion.

For a sublattice $S$ of $L$, we prove the stronger inequality
$\cd(L)\leq2^{-d}\cd(S)$, where $d$ is the larger of the increases in
the numbers of join-reducible and meet-reducible elements. The proof
follows congruence labels along the chains into which prime intervals
of $S$ split in $L$.

We also prove that $\cd(L)=\cd(S)$ if and only if every element of
$L-S$ is comparable with every element of $L$. The theorem of
G. Cz\'edli and L.~Szab\'o relating compatible quasiorders to pairs of
congruences gives the corresponding quasiorder-density results.
\end{abstract}
\maketitle

\section{Introduction}

A chain with $n$ elements has $2^{n-1}$ congruences. Each of its $n-1$
edges may be collapsed or retained, independently. Freese's bound says
that no $n$-element lattice has more congruences, and equality
characterizes chains; see \cite[Lemma~1(3), p.~3458]{Freese1997} and the
short deduction in Section~2.

G\'abor Cz\'edli used this bound to introduce the \emph{congruence
density} of a finite lattice $L$:
\[
   \cd(L)=\frac{|\Con L|}{2^{|L|-1}}.
\]
Here $\Con L$ is the congruence lattice of $L$, and $|L|$ is the number
of its elements. Thus $0<\cd(L)\leq1$. The normalization separates two
questions that raw congruence counts do not distinguish: how large is
the lattice, and how freely can its congruences be chosen?

Cz\'edli determined the second largest congruence count and the
lattices attaining it in \cite{CzedliMany}. Mure\c{s}an and Kulin
determined the next three largest counts and the corresponding lattice
structures in \cite{MuresanKulin}. Cz\'edli's subsequent investigations
\cite{CzedliRanked,CzedliSublatticeCongruences,Czedli332} show that
density carries structural information.

A natural test is sublattice inclusion. If a lattice already contains a
sublattice with few congruences relative to its size, can adjoining
further elements restore the lost density?

Cz\'edli proved that $\cd(L)\leq\cd(S)$ for a sublattice $S$ of a finite
lattice $L$ in three cases: $S$ is cover-preserving; $L$ is modular; or
$L$ can be obtained from $S$ by adjoining one element at a time, obtaining a
sublattice at every step. He asked whether the inequality holds without
these restrictions \cite[Section~1 and Theorem~2.1]{CzedliSublatticeCongruences}.

There is a difficulty even before the counting begins. Restriction sends
a congruence of $L$ to one of $S$, but distinct congruences may have the
same restriction, and some congruences of $S$ may have no extension to
$L$. In particular, we cannot simply multiply $|\Con S|$ by a presumed
number of independent choices for the new elements. We must identify
which choices a congruence actually makes.

Even a one-element chain extension shows why a uniform bound on each
set of extensions is insufficient. Take $S=\{0,1\}$ inside the chain
$L=\{0<c<1\}$. The identity congruence of $S$ has three extensions:
the identity, the congruence collapsing $0,c$, and the congruence
collapsing $c,1$. The universal congruence has one extension. The total
is $3+1=4=2|\Con S|$, although one congruence has more than two
extensions. The desired inequality must balance the whole count.

Our main result does this quantitatively. Let $d_J$ be the number of
join-reducible elements of $L$ minus their number in $S$, and define
$d_M$ dually. Theorem~\ref{prop:monotone-defect} gives
\[
   \cd(L)\leq2^{-\max(d_J,d_M)}\cd(S).
\]
An increase in either reducible-element count therefore has a visible
cost in the bound. These differences are nonnegative, because
reducibility of an element of $S$ is preserved in $L$.

The proof extends Cz\'edli's prime-interval counting method. A prime
interval of $S$ need not remain prime in $L$. We replace it by a maximal
chain and record which edge labels lie below a given congruence.
Repeated labels give equal recorded values; order relations between
labels give implications between values. A counting lemma converts
these restrictions into the factor $2^{-d_J}$. Duality gives the factor
$2^{-d_M}$.

The bound leaves open a second question: when is no density lost at all?
Corollary~\ref{thm:monotonicity} answers it. Equality holds precisely when
every element of $L-S$ is comparable with every element of $L$.
Thus equality permits the insertion of chain elements that introduce no
new incomparability. The proof first treats one deleted element by counting the possible
extensions of each congruence.
The quantitative theorem then reduces the general case to that one.

Section~2 collects the congruence facts used in the argument.
Section~3 proves the quantitative bound, and Section~4 treats equality.
The last section applies the identity $\qd(L)=\cd(L)^2$, due to the
theorem of Cz\'edli and Szab\'o \cite{CzedliSzabo}.

\section{Background}

All lattices and sublattices are finite and nonempty. A sublattice need not contain the
bottom or top of the ambient lattice. We write $0_A$ and $1_A$ for the
bottom and top of a lattice $A$, and $a\prec_A b$ when $b$ covers $a$
in $A$. We omit the subscript when the lattice is clear.

For $a,b\in A$, write $\con_A(a,b)$ for the least congruence of $A$
identifying $a$ and $b$. We order congruences by inclusion. Thus
$\balpha\leq\bbeta$ means that every pair identified by $\balpha$ is
also identified by $\bbeta$.

A \emph{prime interval} is an interval $[a,b]$ with $a\prec b$. Its
\emph{label} is $\con_A(a,b)$. A congruence class, also called a \emph{block}, is a convex
sublattice. Convexity means that $a\equiv b$ and $a\leq x\leq b$ imply
$a\equiv x\equiv b$.

\subsection*{Join-irreducible elements and congruence labels}

Let $\Ji(A)$ be the set of nonzero join-irreducible elements of
$A$, with the inherited order when it is nonempty. Each
$x\in\Ji(A)$ has a unique lower cover, denoted by $x_*$.
An element is join-reducible if it is the join of two strictly smaller
elements, or equivalently, if it has at least two lower covers.
Meet-reducibility is defined dually. Every element other than $0_A$ is
either join-irreducible or join-reducible.

For a finite nonempty ordered set $P$, let $\Id(P)$ be the lattice of
its order ideals under inclusion. A subset $I\subseteq P$ is an order ideal if
$q\in I$ and $p\leq q$ imply $p\in I$. In particular,
$\varnothing\in\Id(P)$.

The congruence lattice $\Con A$ is distributive. If $|A|\geq2$, then
$\Ji(\Con A)$ is nonempty, and its usual ideal representation is
\begin{equation}\label{eq:ideal-representation}
   \Con A\cong\Id(\Ji(\Con A)),\qquad
   \btheta\longmapsto\{p\in\Ji(\Con A):p\leq\btheta\}.
\end{equation}
The one-element lattice has one congruence and will be treated directly
when needed.

The labels $\con_A(x_*,x)$, for $x\in\Ji(A)$, exhaust
$\Ji(\Con A)$; different $x$ may give the same label
\cite[Lemma~1]{Freese1997}; see also \cite[Chapter~III]{GG11}.
Consequently, for each congruence $\btheta$, the function on $\Ji(A)$
whose value at $x$ is $1$ when $\con_A(x_*,x)\leq\btheta$ and $0$
otherwise determines $\btheta$. We shall use such functions in
Section~3. Not every function to $\{0,1\}$ need arise from a congruence:
equal labels force equal values, and order relations between labels
force implications between values.

We use two standard consequences repeatedly. First, every prime-interval
label is join-irreducible in $\Con A$. Second, join-irreducible elements
of a finite distributive lattice are join-prime: if
$p\leq q_1\vee\cdots\vee q_m$, then $p\leq q_i$ for some $i$.
For example, distributivity gives
$p=(p\wedge q_1)\vee\cdots\vee(p\wedge q_m)$, and
join-irreducibility then gives the assertion.

If $a=c_0\prec c_1\prec\cdots\prec c_m=b$ is a maximal chain in an
interval, then
\begin{equation}\label{eq:chain-labels}
   \con_A(a,b)=\bigvee_{i=1}^m\con_A(c_{i-1},c_i).
\end{equation}
Indeed, collapsing the endpoints collapses every edge by convexity.
Conversely, collapsing all edges identifies the endpoints by transitivity.
This identity explains why chains can replace prime intervals when we
pass to a larger lattice.

\subsection*{The elementary bound and a transposition}

For $|A|\geq2$, equation~\eqref{eq:ideal-representation} and the
surjection of labels give
\[
   |\Con A|\leq2^{|\Ji(\Con A)|}
      \leq2^{|\Ji(A)|}\leq2^{|A|-1}.
\]
These inequalities also hold for the one-element lattice, since every
term is then $1$.
If $A$ is not a chain, choose incomparable elements $u,v$. Their join is
join-reducible, so $|\Ji(A)|\leq|A|-2$ and the last bound is strict.
In a chain, any selection of collapsed edges defines a congruence.
This proves both the bound and its equality statement.

The following elementary transposition replaces an arbitrary prime
interval by one with a join-irreducible upper endpoint and the same label.
It is the argument used by Day \cite[p.~71]{Day1979}; see also
\cite[Lemma~3.4]{CzedliSublatticeCongruences}.

\begin{lemma}\label{lem:edge-transposition}
Let $a\prec b$ in a finite lattice $A$. Choose $f$ minimal subject to
$f\leq b$ and $f\not\leq a$. Then $f\in\Ji(A)$,
$f_*=f\wedge a$, and $f\vee a=b$. The prime intervals $[f_*,f]$ and
$[a,b]$ have the same label.
\end{lemma}

\begin{proof}
Such an $f$ exists because $b\not\leq a$. If $z<f$, minimality gives
$z\leq a$, hence $z\leq f\wedge a$. Since $f\not\leq a$, the element
$f\wedge a$ is strictly below $f$. It is therefore the greatest element
below $f$, proving that $f$ is join-irreducible and $f_*=f\wedge a$.
Also $a<f\vee a\leq b$, so $f\vee a=b$.

If a congruence identifies $f_*$ with $f$, joining with $a$ identifies
$a=f_*\vee a$ with $b=f\vee a$. Conversely, if it identifies $a$
with $b$, meeting with $f$ identifies $f_*=a\wedge f$ with $f=b\wedge f$.
Thus the two pairs generate the same congruence.
\end{proof}

\subsection*{Glued sums}

The \emph{glued sum} of $A$ and $B$ identifies $1_A$ with $0_B$ and
places $A$ below $B$. It is the sublattice
\[
   (A\times\{0_B\})\cup(\{1_A\}\times B)
\]
of $A\times B$. Any pair of factor congruences extends by restricting
their product congruence to this sublattice. The extension is unique:
if a class meets both factors, convexity forces it to contain the common
element, so its parts are determined by the two restrictions.

Congruence counts therefore multiply under gluing. Since the glued sum
has $|A|+|B|-1$ elements, densities multiply as well. In particular,
gluing a chain at either end leaves density unchanged. We will use the
same observation when a chain segment occurs between two lattice pieces.

\section{Density under sublattice inclusion}\label{sec:monotonicity}

Let $S\subseteq L$ be a sublattice inclusion. Write $d_J$ for the number
of join-reducible elements of $L$ minus their number in $S$, and define
$d_M$ dually. If $s=u\vee v$ witnesses join-reducibility in $S$, the same
equality holds in $L$. Thus the join-reducible elements of $S$ form a
subset of those of $L$, and $d_J\geq0$. Similarly, $d_M\geq0$.

\begin{theorem}\label{prop:monotone-defect}
For every sublattice $S$ of a finite lattice $L$,
\[
   \cd(L)\leq2^{-\max(d_J,d_M)}\cd(S).
\]
\end{theorem}

If $r=|L-S|$, the estimate with $d_J$ is equivalent to
\begin{equation}\label{eq:count-target}
   |\Con L|\leq2^{r-d_J}|\Con S|.
\end{equation}
The exponent has a useful interpretation:
\[
   r-d_J=|\Ji(L)|-|\Ji(S)|.
\]
It measures the change in the number of join-irreducible elements.
The proof must show that the order relations between congruence labels
in $S$ still impose enough restrictions on congruences of $L$.

The factor can be attained. A three-element chain $S$ inside the
four-element Boolean lattice $L$ has $d_J=d_M=1$, and
$\cd(L)=1/2=2^{-1}\cd(S)$. Thus the additional reducibility can account
for the entire decrease.

The reducible-element counts alone do not determine equality. Let $B_4$
be the four-element Boolean lattice and $N_5$ the pentagon. Deleting an
interior element of the longer arm of $N_5$ leaves $B_4$. Both lattices
have one join-reducible and one meet-reducible element, but
\[
   \cd(N_5)=\frac5{16}<\frac12=\cd(B_4).
\]
For clarity, write the pentagon as $0<a<b<1$ and $0<c<1$, with $c$
incomparable with $a,b$. Its five congruences have nontrivial blocks
$\{a,b\}$, or the two blocks $\{0,a,b\},\{c,1\}$, or the two blocks
$\{0,c\},\{a,b,1\}$, in addition to the identity and universal
congruences. This also checks the stated density directly.

\subsection*{A counting lemma}

We first specify the objects to be counted. If $B$ is a finite set of
indices, a \emph{binary assignment on $B$} is a function
\[
   v\colon B\longrightarrow\{0,1\}.
\]
The value $v(i)$ is the coordinate at the index $i\in B$. An implication
from an index $j$ to an index $i$ means the requirement
\[
   v(j)=1\quad\Longrightarrow\quad v(i)=1,
   \qquad\text{equivalently,}\qquad v(j)\leq v(i).
\]
The order on the values is $0<1$.

We partition $B$ into nonempty sets $B_p$, one for each point $p$ of a
finite ordered set $P$. Thus the $B_p$ are sets of indices. They are not
sets of binary values. Pairwise disjointness means
$B_p\cap B_q=\varnothing$ whenever $p\ne q$; it places no restriction
on the values assigned to indices in different sets.

\smallskip
\noindent\emph{Example.}
Let $P=\{p,q\}$ with $p<q$, and take
\[
   B_p=\{1,2\},\qquad B_q=\{3\},\qquad B=\{1,2,3\}.
\]
Require both implications $v(3)\leq v(1)$ and $v(3)\leq v(2)$.
Writing an assignment as the triple $(v(1),v(2),v(3))$, the possible
assignments are exactly
\[
   (0,0,0),\quad(1,0,0),\quad(0,1,0),\quad(1,1,0),\quad(1,1,1).
\]
If we also require $v$ to be constant on $B_p$, only
$(0,0,0)$, $(1,1,0)$, and $(1,1,1)$ remain. They correspond to the
three order ideals $\varnothing$, $\{p\}$, and $\{p,q\}$ of $P$.
The next lemma bounds how much is lost by imposing constancy on every
$B_p$. In this example its bound is $2^{3-2}\cdot3=6$, while the
actual count is five.

\begin{lemma}\label{lem:monotone-groups}
Let $P$ be a finite nonempty ordered set. For each $p\in P$, let $B_p$
be a nonempty finite set, and suppose that these sets are pairwise
disjoint. Put
\[
   B=\bigcup_{p\in P}B_p,
   \qquad M=|B|=\sum_{p\in P}|B_p|.
\]
Let $R$ be a relation on $B$ such that
\begin{equation}\label{eq:implication-relation}
   R\subseteq\bigcup_{p<q}(B_q\times B_p).
\end{equation}
Assume that, for every $p<q$ and every $i\in B_p$, there is a
$j\in B_q$ with $(j,i)\in R$.

Call a function $v\colon B\to\{0,1\}$ \emph{admissible} if
$v(j)\leq v(i)$ for every $(j,i)\in R$. Then the number of admissible
functions is at most
\[
   2^{M-|P|}|\Id(P)|.
\]
\end{lemma}

\begin{proof}
List the points of $P$ as $p_1,\ldots,p_k$, where $k=|P|$, so that
every point occurs before all points strictly below it. Such a list is
obtained by successively removing a maximal point.

Let $\mathcal A_0$ be the set of all admissible functions. For
$s=1,\ldots,k$, let $\mathcal A_s$ consist of those functions in
$\mathcal A_{s-1}$ that are constant on $B_{p_s}$. We prove that
\begin{equation}\label{eq:one-constancy-step}
   |\mathcal A_{s-1}|
      \leq 2^{|B_{p_s}|-1}|\mathcal A_s|.
\end{equation}

Fix $s$, put $p=p_s$ and $b=|B_p|$, and fix values on $B-B_p$ that
extend to at least one function in $\mathcal A_{s-1}$. We compare the
number of extensions before and after requiring constancy on $B_p$.
Every $B_q$ with $q>p$ already has a constant value, by our choice of
the list.

If some such $B_q$ has value $1$, the hypothesis on $R$ forces
$v(i)=1$ for every $i\in B_p$. There is exactly one extension, and it
is constant on $B_p$.

Suppose instead that every $B_q$ with $q>p$ has value $0$. No
implication entering $B_p$ then forces a value. An implication leaving
$B_p$ can only force its source to have value $0$, when its target has
the fixed value $0$. There are no implications within $B_p$, by
\eqref{eq:implication-relation}.

If no index in $B_p$ is forced to have value $0$, there are $2^b$
extensions before the new requirement and two afterward. If at least
one is forced to have value $0$, there are at most $2^{b-1}$ extensions
before the requirement, and exactly one constant extension, the
all-zero function on $B_p$, afterward. Thus in every case the first
count is at most $2^{b-1}$ times the second. Summing over the fixed
functions on $B-B_p$ proves \eqref{eq:one-constancy-step}.

A function in $\mathcal A_k$ is constant on each $B_p$. Let $I$ be
the set of points $p$ whose common value is $1$. If $p<q$ and $q\in I$,
the hypothesis on $R$ forces the value $1$ throughout $B_p$, so
$p\in I$. Hence $I$ is an order ideal. Conversely, an order ideal $I$
gives an admissible function by assigning value $1$ precisely to the
indices in $\bigcup_{p\in I}B_p$. These constructions are inverse,
and therefore $|\mathcal A_k|=|\Id(P)|$.

Multiplying \eqref{eq:one-constancy-step} for $s=1,\ldots,k$ gives
\[
   |\mathcal A_0|
      \leq\prod_{p\in P}2^{|B_p|-1}|\mathcal A_k|
      =2^{M-|P|}|\Id(P)|.
\]
\end{proof}

\subsection*{Proof of the quantitative inequality}

In the application, the original index set will be $\Ji(L)$.
Different indices may carry the same congruence label; the corresponding
values then agree for every congruence. We will delete specified
duplicates and define the sets $B_p$ as subsets of the remaining index
set. This keeps the indices, their labels, and their assigned values
separate throughout the argument.

\begin{proof}[Proof of Theorem~\ref{prop:monotone-defect}]
It suffices to prove \eqref{eq:count-target}, the estimate using $d_J$.
Applying the same argument to the dual lattices gives the estimate
using $d_M$, and we take the stronger of the two bounds. Put
$r=|L-S|$.

If $|S|=1$, then $|\Con S|=1$ and $|\Ji(S)|=0$. The elementary
bound in Section~2 gives
\[
   |\Con L|\leq2^{|\Ji(L)|}=2^{r-d_J}|\Con S|.
\]
We may therefore assume that $|S|\geq2$. The ordered sets
\[
   P=\Ji(\Con S),\qquad Q=\Ji(\Con L)
\]
are then nonempty. We organize the proof in five steps.

\smallskip
\noindent\emph{1. Transfer the labels of $S$ to $L$.}
For each $x\in\Ji(S)$, let $a_x$ be its lower cover in $S$, and choose
a maximal chain
\[
   a_x=c_{x,0}\prec_L c_{x,1}\prec_L\cdots
       \prec_L c_{x,m_x}=x
\]
in $L$. Put $p_x=\con_S(a_x,x)$. The labels $p_x$ exhaust $P$.

Suppose that $p_x\leq p_y$. The restriction of
$\con_L(a_y,y)$ to $S$ is a congruence containing $p_y$, and hence
$p_x$. It identifies $a_x$ with $x$, so
\[
   \con_L(a_x,x)\leq\con_L(a_y,y).
\]
In particular, equal labels in $S$ give equal principal congruences in
$L$. Thus the definition
\[
   \blambda_p=\con_L(a_x,x)\quad\text{when }p_x=p
\]
does not depend on the choice of $x$. The map
$p\mapsto\blambda_p$ is isotone, but need not be injective.
Every $\blambda_p$ is nonzero because $a_x\ne x$.

\smallskip
\noindent\emph{2. Represent congruences by functions on $\Ji(L)$.}
Define the label map
\[
   \ell\colon\Ji(L)\longrightarrow Q,
   \qquad \ell(f)=\con_L(f_*,f).
\]
For $\btheta\in\Con L$, define the binary assignment
\[
   e_{\btheta}\colon\Ji(L)\longrightarrow\{0,1\},
   \qquad
   e_{\btheta}(f)=
   \begin{cases}
      1,&\ell(f)\leq\btheta,\\
      0,&\ell(f)\not\leq\btheta.
   \end{cases}
\]
The label map is onto $Q$, so
\eqref{eq:ideal-representation} shows that $e_{\btheta}$ determines
$\btheta$. Thus $\btheta\mapsto e_{\btheta}$ is injective.

For every edge $[c_{x,i-1},c_{x,i}]$ of every chosen chain, choose
$f_{x,i}$ minimal subject to
\[
   f_{x,i}\leq c_{x,i},\qquad f_{x,i}\not\leq c_{x,i-1}.
\]
Lemma~\ref{lem:edge-transposition} gives $f_{x,i}\in\Ji(L)$ and
\begin{equation}\label{eq:edge-index-label}
   \ell(f_{x,i})=\con_L(c_{x,i-1},c_{x,i}).
\end{equation}
We now check that all the indices $f_{x,i}$ are distinct.

For two edges of the same chain, with $i<j$, we have
$f_{x,i}\leq c_{x,i}\leq c_{x,j-1}$, whereas
$f_{x,j}\not\leq c_{x,j-1}$. Hence $f_{x,i}\ne f_{x,j}$.

For distinct $x,y\in\Ji(S)$, suppose that an index $f$ belonged to
both chains. Then $f\leq x\wedge y$. The meet is in $S$ and is
strictly below at least one of $x,y$, say $x$. Since $a_x$ is the
greatest element of $S$ strictly below $x$, we have $x\wedge y\leq a_x$.
But an index selected for the chain from $a_x$ to $x$ cannot lie below
$a_x$. This contradiction proves the claim.

\smallskip
\noindent\emph{3. Delete specified duplicate indices.}
For $p\in P$, define
\[
   X_p=\{x\in\Ji(S):p_x=p\},\qquad t_p=|X_p|.
\]
The sets $X_p$ form a partition of $\Ji(S)$ into nonempty sets.
Also let
\[
   T_p=\max\{u\in Q:u\leq\blambda_p\},\qquad b_p=|T_p|,
\]
where $\max$ denotes the set of maximal elements. Since
$\blambda_p\ne0$, this finite set is nonempty.

We claim that every $u\in T_p$ occurs as an edge label on every chain
whose top $x$ belongs to $X_p$. By \eqref{eq:chain-labels}, the edge
labels on such a chain have join $\blambda_p$. Join-primality gives an
edge label $u'$ with $u\leq u'$. Since $u'\in Q$ and
$u'\leq\blambda_p$, maximality of $u$ gives $u=u'$, proving the
claim.

For each $x\in X_p$ and $u\in T_p$, choose one edge with label $u$.
Write $g_{x,u}$ for the index $f_{x,i}$ already selected for that edge
in Step~2. Then $\ell(g_{x,u})=u$. All selected indices $g_{x,u}$ are
distinct: Step~2 distinguishes edges on different chains, and different
labels cannot be carried by the same edge.

Choose one $x_p\in X_p$ for each $p$. Define
\begin{align*}
   B_p&=\{g_{x_p,u}:u\in T_p\},\\
   D_p&=\{g_{x,u}:x\in X_p-\{x_p\},\ u\in T_p\},\\
   C&=\Ji(L)-\bigcup_{p\in P}D_p.
\end{align*}
The set $C$ is the retained index set. The sets $B_p$ are pairwise
disjoint subsets of $C$, and $|B_p|=b_p$. The set $D_p$ consists of
the $(t_p-1)b_p$ indices that we delete for $p$.

For every congruence $\btheta$ and every deleted index $g_{x,u}$,
\[
   e_{\btheta}(g_{x,u})=e_{\btheta}(g_{x_p,u}),
\]
because both indices have label $u$. Thus $e_{\btheta}|_C$ recovers
the values at all deleted indices and therefore recovers
$e_{\btheta}$. By Step~2, the map
$\btheta\mapsto e_{\btheta}|_C$ remains injective.

For example, if $t_p=3$ and $b_p=2$, we have selected six distinct
indices. We retain the two from the chain with top $x_p$ and delete
the four from the other two chains. Each deleted value is recovered
from the retained index with the same label.

Put $N=|C|$. Since the deleted indices are all distinct, we obtain
\begin{align}
   N&=|\Ji(L)|-\sum_{p\in P}(t_p-1)b_p \notag\\
    &\leq|\Ji(L)|-\sum_{p\in P}(t_p-1) \notag\\
    &=|\Ji(L)|-|\Ji(S)|+|P| \notag\\
    &=r-d_J+|P|.\label{eq:number-coordinates}
\end{align}
Here $b_p\geq1$ and $\sum_p t_p=|\Ji(S)|$. For the last equality,
subtract the identities
\[
   |\Ji(A)|=|A|-1-\#\{\text{join-reducible elements of }A\}
\]
for $A=L$ and $A=S$. Each identity excludes the bottom of its own
lattice, so the subtraction is valid even when $0_S\ne0_L$.

\smallskip
\noindent\emph{4. Specify the implication relation.}
Put $B=\bigcup_{p\in P}B_p$ and $M=|B|$. If $p<q$ and $i\in B_p$,
then
\[
   \ell(i)\leq\blambda_p\leq\blambda_q.
\]
In the finite ordered set $\{u\in Q:u\leq\blambda_q\}$, choose a
maximal element $u$ above $\ell(i)$. Then $u\in T_q$, and
$j=g_{x_q,u}\in B_q$ satisfies $\ell(i)\leq\ell(j)$.
Consequently, every congruence assignment satisfies
\[
   e_{\btheta}(j)\leq e_{\btheta}(i).
\]

For each $p<q$ and each $i\in B_p$, make one such choice of $j$ and
put $(j,i)$ in $R$. These pairs define a relation satisfying every
hypothesis of Lemma~\ref{lem:monotone-groups}. In particular,
$R\subseteq\bigcup_{p<q}(B_q\times B_p)$, and each required target
index receives an implication from the specified upper set.

We now count all functions on $C$ that satisfy the implications in
$R$. Every function $e_{\btheta}|_C$ is among them. Other restrictions
may be ignored for an upper bound. In particular, indices in different
$B_p$ may have equal labels, but we need not impose all equalities
between their assigned values. The indices themselves remain distinct.

\smallskip
\noindent\emph{5. Count and normalize.}
Lemma~\ref{lem:monotone-groups} bounds the admissible functions on $B$
by $2^{M-|P|}|\Id(P)|$. Values on $C-B$ are unrestricted in our upper
bound, contributing a factor $2^{N-M}$. The injection from Step~3
therefore gives
\[
   |\Con L|\leq2^{N-M}2^{M-|P|}|\Id(P)|
      =2^{N-|P|}|\Id(P)|
      \leq2^{r-d_J}|\Con S|.
\]
The last inequality uses \eqref{eq:number-coordinates} and the ideal
representation of $\Con S$. This proves \eqref{eq:count-target}.
Dividing by $2^{|L|-1}$ gives
$\cd(L)\leq2^{-d_J}\cd(S)$, as required.
\end{proof}

\section{Equality of densities}\label{sec:equality}

The quantitative theorem shows that equality of densities forces
$d_J=d_M=0$. The pentagon example shows that the converse fails.
We need to determine what else can make the inequality strict.
We begin with a sublattice obtained by deleting one element.

For $\btheta\in\Con L$, its \emph{restriction} to $S$ is the relation
$\btheta|_S=\btheta\cap(S\times S)$, which belongs to $\Con S$.
We say that $\btheta$ \emph{extends} $\balpha\in\Con S$ if
$\btheta|_S=\balpha$. An extension need not exist.

\begin{lemma}\label{lem:monotone-one}
Let $S=L-\{c\}$ be a nonempty sublattice of a finite lattice $L$.
Then $\cd(L)=\cd(S)$ if and only if $c$ is comparable with every
element of $L$.

If $c$ is neither $0_L$ nor $1_L$, then it has a unique lower cover
$a$ and a unique upper cover $b$. If there is an $x\in S$ with
$a<x<b$, then
\[
   \cd(L)\leq\tfrac12\cd(S).
\]
\end{lemma}

\begin{proof}
Since $|L|=|S|+1$, equality of densities is equivalent to
$|\Con L|=2|\Con S|$.

If $c=0_L$, the least element $0_S$ of the nonempty sublattice $S$
lies below every element of $L$ other than $c$. Thus $L$ is obtained
by gluing the two-element chain $\{c,0_S\}$ below $S$.
Its number of congruences is $2|\Con S|$.
The case $c=1_L$ is dual. In both cases, $c$ is comparable with
every element of $L$, as required.

Assume that $c$ is neither endpoint. It has a unique lower cover $a$.
Indeed, two distinct lower covers would belong to $S$ and have join
$c$, contrary to closure of $S$ under joins. Dually, $c$ has a unique
upper cover $b$. Every element strictly below $c$ is therefore at most
$a$, and every element strictly above $c$ is at least $b$.

\smallskip
\noindent\emph{Equality under the comparability condition.}
Suppose that $c$ is comparable with every element of $L$.
Then every element of $S$ lies below $a$ or above $b$.
Put $A=[0_L,a]$ and $B=[b,1_L]$.
The lattice $L$ is obtained by gluing $A$, the chain $a<c<b$, and
$B$, in this order. The lattice $S$ is obtained by gluing the same
two intervals with the chain $a<b$ between them.
Consequently,
\[
   |\Con L|=4|\Con A||\Con B|,
   \qquad
   |\Con S|=2|\Con A||\Con B|.
\]
This proves equality of densities in this case.

For the converse, assume from now on that $c$ has an incomparable
element. We shall prove $|\Con L|<2|\Con S|$.

\smallskip
\noindent\emph{The possible extensions of a fixed congruence.}
For $\balpha\in\Con S$, let $e(\balpha)$ be the number of its
extensions to $L$. Since every congruence of $L$ restricts to exactly
one congruence of $S$,
\begin{equation}\label{eq:extension-count}
   |\Con L|=\sum_{\balpha\in\Con S}e(\balpha).
\end{equation}

If $\balpha$ identifies $a$ with $b$, convexity forces $c$ into their
common block in every extension. All other blocks are already fixed
by $\balpha$, so $e(\balpha)\leq1$.

If $\balpha$ separates $a$ and $b$, an extension has at most three
choices for the block of $c$: a singleton, the block containing $a$,
or the block containing $b$.
To justify this list, suppose that an extension identifies $c$ with
some $z\in S$. If $z<c$, convexity identifies $c$ with $a$; if
$z>c$, it identifies $c$ with $b$.
If $z$ is incomparable with $c$, compatibility with meet and join
identifies $c$ with both $c\wedge z<c$ and $c\vee z>c$.
Convexity then identifies $a$ with $b$, contrary to the assumption.
Once the block of $c$ is chosen, the extending equivalence relation
is determined. Some of these choices may fail to be congruences;
we only use them as an upper bound.

There are now three cases. Either the interval $[a,b]$ in $S$ has an
interior element, or $a\prec_S b$. In the latter case, we distinguish
whether this cover is unique at both endpoints.

\smallskip
\noindent\emph{Case 1: there is an $x\in S$ with $a<x<b$.}
The element $x$ is incomparable with $c$, since $a\prec_L c\prec_L b$.
The unique-cover properties of $c$ give
$c\wedge x=a$ and $c\vee x=b$.
For any extension $\btheta$ of $\balpha$, compatibility gives
\[
   c\equiv a\pmod{\btheta}
      \quad\Longleftrightarrow\quad x\equiv b\pmod{\balpha},
\]
and, dually,
\[
   c\equiv b\pmod{\btheta}
      \quad\Longleftrightarrow\quad a\equiv x\pmod{\balpha}.
\]
For the first equivalence, join $c\equiv a$ with $x$ in the forward
direction, and meet $x\equiv b$ with $c$ in the reverse direction.

If $\balpha$ separates $a$ and $b$, the two right-hand conditions
cannot both hold. If one holds, it determines which endpoint block
contains $c$. If neither holds, $c$ must be a singleton.
Thus $e(\balpha)\leq1$ for every $\balpha$, including those that
identify $a$ and $b$. Equation~\eqref{eq:extension-count} gives
\[
   |\Con L|\leq|\Con S|,
   \qquad \cd(L)\leq\tfrac12\cd(S).
\]
This also proves the additional assertion of the lemma.

\smallskip
\noindent\emph{Case 2: $a\prec_S b$, and an endpoint has another cover.}
Suppose first that $a$ has an upper cover $d\ne b$ in $S$.
The two covers satisfy $b\wedge d=a$.
We also have $c\vee d=b\vee d$.
Indeed, $d\not\leq c$, since otherwise $d\leq a$.
Hence $c\vee d>c$, so $b\leq c\vee d$ by the unique upper-cover
property of $c$. This gives $b\vee d\leq c\vee d$; the reverse
inequality follows from $c\leq b$.

If an extension identifies $c$ with $a$, joining with $d$ gives
$b\vee d\equiv d$. Meeting this pair with $b$ gives $b\equiv a$.
Therefore, when $\balpha$ separates $a$ and $b$, the block of $c$
cannot be the block of $a$. There are at most two extensions.
Dually, if $b$ has another lower cover in $S$, the block of $c$
cannot be the block of $b$, and again there are at most two extensions.

Congruences identifying $a$ and $b$ have at most one extension.
There is at least one such congruence, the universal congruence of $S$.
Consequently,
\[
   |\Con L|\leq2|\Con S|-1<2|\Con S|.
\]
The density inequality is strict.

\smallskip
\noindent\emph{Case 3: $b$ is the unique upper cover of $a$ in $S$,
and $a$ is the unique lower cover of $b$ in $S$.}
The preliminary bound allows three extensions when $a$ and $b$ are
separated, and only one when they are identified.
We shall show that more congruences identify this pair than separate it.

Every element of $S$ lies below $a$, lies above $b$, or is incomparable
with both. For example, an element strictly above $a$ must be at least
its unique upper cover $b$; the statement below $b$ is dual.
If $y$ is incomparable with both endpoints, these same properties give
\begin{equation}\label{eq:equal-meets-joins}
   a\vee y=b\vee y,\qquad a\wedge y=b\wedge y.
\end{equation}
For the first equality, $a\vee y>a$, so $b\leq a\vee y$.
It follows that $b\vee y\leq a\vee y$; the reverse inequality follows
from $a\leq b$. The second equality is dual.

Let $\bbeta$ be the equivalence relation on $S$ whose only
nonsingleton block is $\{a,b\}$.
It is a congruence. To check this, it suffices to check compatibility
of the pair $a,b$ with each $y\in S$.
For $y\leq a$, the meets are equal and the joins are $a,b$;
for $y\geq b$, the joins are equal and the meets are $a,b$.
For an element incomparable with both, use
\eqref{eq:equal-meets-joins}.
Thus $\bbeta=\con_S(a,b)$.
Moreover, $\bbeta$ is an atom of $\Con S$, since the identity relation
is its only proper equivalence refinement.

Put $R=\Ji(\Con S)$ and define
\[
   \mathcal I_0=\{I\in\Id(R):\bbeta\notin I\},
   \qquad
   \mathcal I_1=\{I\in\Id(R):\bbeta\in I\}.
\]
Let $N_0=|\mathcal I_0|$ and $N_1=|\mathcal I_1|$.
Under the representation~\eqref{eq:ideal-representation}, these are
the numbers of congruences of $S$ that separate $a,b$ and identify
$a,b$, respectively. Hence
\begin{equation}\label{eq:atom-extension-bound}
   |\Con S|=N_0+N_1,
   \qquad |\Con L|\leq3N_0+N_1.
\end{equation}

Since $\bbeta$ is an atom, it is a minimal element of the ordered
set $R$. Thus adding $\bbeta$ to an ideal that omits it still gives
an ideal. The map
\[
   \Phi\colon\mathcal I_0\longrightarrow\mathcal I_1,
   \qquad \Phi(I)=I\cup\{\bbeta\},
\]
is injective: deleting $\bbeta$ recovers $I$. This proves $N_0\leq N_1$.
We shall exhibit an ideal in $\mathcal I_1$ outside its image.

Choose $y\in S$ incomparable with $c$.
The partition of $S$ described above shows that $y$ is incomparable
with both $a$ and $b$. Put
\[
   \balpha=\con_S(a\wedge y,y).
\]
Joining the generating pair with $a$ gives
$a\equiv a\vee y\pmod{\balpha}$.
Joining it with $b$ gives $b\equiv b\vee y\pmod{\balpha}$.
The right-hand elements are equal by~\eqref{eq:equal-meets-joins}.
Hence $a\equiv b\pmod{\balpha}$, or $\bbeta\leq\balpha$.

By~\eqref{eq:chain-labels}, $\balpha$ is the join of the labels on a
maximal chain from $a\wedge y$ to $y$ in $S$.
The atom $\bbeta$ is join-prime in the distributive lattice $\Con S$,
so one of these labels, say $\btheta\in R$, satisfies
$\bbeta\leq\btheta$.
In fact, $\bbeta<\btheta$.
Otherwise $\bbeta$ would collapse an edge of that chain.
Its only nonsingleton block is $\{a,b\}$, whereas every element of
the chain is at most $y$ and $b\not\leq y$.

Now consider the principal order ideal
\[
   {\downarrow}\btheta
      =\{\blambda\in R:\blambda\leq\btheta\}.
\]
It contains $\bbeta$, so it belongs to $\mathcal I_1$.
It cannot equal $\Phi(I)$ for any $I\in\mathcal I_0$:
such an $I$ would contain $\btheta$ but omit the smaller element
$\bbeta$, and hence would not be an order ideal.
Thus $\Phi$ is not onto, and $N_0<N_1$.

Equation~\eqref{eq:atom-extension-bound} now gives
\[
   |\Con L|\leq3N_0+N_1
      <2(N_0+N_1)=2|\Con S|.
\]
This completes the strict inequality in all three cases, and the proof.
\end{proof}

\begin{corollary}\label{thm:monotonicity}
Let $S$ be a sublattice of a finite lattice $L$. Then
$\cd(L)\leq\cd(S)$, with equality if and only if every element of
$L-S$ is comparable with every element of $L$.
\end{corollary}

\begin{proof}
The inequality is Theorem~\ref{prop:monotone-defect} with
$d_J,d_M\geq0$. If $S=L$, the equality criterion is immediate.

Suppose that $\cd(L)=\cd(S)$ and $S\ne L$.
Since $\cd(S)>0$, the quantitative theorem forces $d_J=d_M=0$.
Every join-reducible element of $S$ remains join-reducible in $L$.
The two sets have the same size, so no element of $L-S$ is
join-reducible in $L$. Dually, none is meet-reducible.

Fix $c\in L-S$. The set $L-\{c\}$ is a sublattice.
Indeed, if the join of two remaining elements were $c$, both would
be strictly below $c$, making $c$ join-reducible.
The meet argument is dual. This sublattice is nonempty because it
contains $S$.

Apply the quantitative theorem to the two inclusions
$S\subseteq L-\{c\}\subseteq L$. It gives
\[
   \cd(L)\leq\cd(L-\{c\})\leq\cd(S)=\cd(L).
\]
Thus deleting $c$ preserves density.
Lemma~\ref{lem:monotone-one} shows that $c$ is comparable with every
element of $L$. Since $c$ was arbitrary, this proves necessity.

Conversely, suppose that every element of $L-S$ is comparable with
every element of $L$.
Any incomparable pair in $L$ must then lie in $S$, so its join and
meet also lie in $S$.
It follows that the elements of $L-S$ can be deleted in any order
while preserving the sublattice property.
To see this at one step, consider two remaining elements.
If they are comparable, their join and meet are those two elements.
If they are incomparable, both belong to $S$, as do their join and
meet. In either case, the element being deleted is neither their join
nor their meet.

Each element being deleted is still comparable with every element of
the current lattice. Lemma~\ref{lem:monotone-one} therefore preserves
density at each step. The last lattice is $S$, proving
$\cd(L)=\cd(S)$.
\end{proof}

The criterion distinguishes the examples from Section~\ref{sec:monotonicity}.
Adjoining a new bottom or top preserves density, as does inserting a
chain segment between two pieces of a glued sum.
Subdividing an arbitrary edge need not preserve density: the new
element may be incomparable with elements elsewhere in the lattice.
This is what happens when $B_4$ is extended to $N_5$.

\section{Quasiorder density}

A \emph{quasiorder} on a set $A$ is a reflexive, transitive binary
relation $\rho\subseteq A\times A$.
For a lattice $A$, the relation is \emph{compatible} if, whenever
$(x_1,y_1),(x_2,y_2)\in\rho$, it also contains both pairs
\[
   (x_1\wedge x_2,\,y_1\wedge y_2),
   \qquad
   (x_1\vee x_2,\,y_1\vee y_2).
\]
Write $\Quo A$ for the lattice of compatible quasiorders of $A$,
ordered by inclusion. A congruence is a symmetric compatible
quasiorder; symmetry is not required for an arbitrary member of
$\Quo A$.

Cz\'edli and Szab\'o proved that
\[
   \Quo A\cong\Con A\times\Con A
\]
as lattices~\cite[Theorem~1]{CzedliSzabo}.
For finite $A$, this gives $|\Quo A|=|\Con A|^2$.
The maximum number of compatible quasiorders of an $n$-element
lattice is therefore $2^{2n-2}$, attained by a chain.

Following Cz\'edli~\cite{CzedliSublatticeCongruences}, define the
\emph{quasiorder density} of $A$ by
\[
   \qd(A)=\frac{|\Quo A|}{2^{2|A|-2}}.
\]
The counting identity becomes
\begin{equation}\label{eq:quasiorder-square}
   \qd(A)
      =\left(\frac{|\Con A|}{2^{|A|-1}}\right)^2
      =\cd(A)^2.
\end{equation}
We can therefore obtain both the inequality and its equality criterion
directly from the congruence-density results.

\begin{corollary}\label{cor:quasiorder}
Let $S$ be a sublattice of a finite lattice $L$, and let $d_J,d_M$ be
as in Section~\ref{sec:monotonicity}. Then
\[
   \qd(L)\leq2^{-2\max(d_J,d_M)}\qd(S)\leq\qd(S).
\]
Equality $\qd(L)=\qd(S)$ holds if and only if every element of
$L-S$ is comparable with every element of $L$.
\end{corollary}

\begin{proof}
Square the inequality in Theorem~\ref{prop:monotone-defect} and use
\eqref{eq:quasiorder-square}.
Since congruence densities are positive, equality of their squares is
equivalent to equality of the densities.
Corollary~\ref{thm:monotonicity} gives the stated criterion.
\end{proof}

\end{document}